\documentclass[reqno]{amsart}
\usepackage[left=1in, right=1in, top=1in, bottom=1in]{geometry}

\usepackage[T2A]{fontenc}
\usepackage[utf8]{inputenc}
\usepackage[russian,english]{babel}

\usepackage{amssymb}
\usepackage{amsmath}
\usepackage{amsfonts}
\usepackage{mathrsfs}
\usepackage{amsthm}
\usepackage{hyperref}
\usepackage{graphicx}
\usepackage{float}
\usepackage{enumerate} 
\usepackage{setspace}
\graphicspath{ {./figure/} }

\hypersetup{
  colorlinks,
  linkcolor={red!70!black},
  citecolor={blue!70!black},
  urlcolor={blue!90!black}
  }
\newcommand{\N}{\ensuremath{\mathbb{N}}}

\newcommand{\R}{\ensuremath{\mathbb{R}}}
\newcommand{\h}{\ensuremath{\mathbb{H}}}

\newcommand{\Div}{\ensuremath{\text{div}}}
\newcommand{\Exc}{\ensuremath{\text{Exc}}}

\newcommand{\gr}{\ensuremath{\text{gr}}}

\newcommand{\symdiff}{\mathop{\triangle}}

\newtheorem{theorem}{Theorem}[section]
\newtheorem*{theorema}{Theorem A}

\newtheorem*{remark}{Remark}

\theoremstyle{definition}
\newtheorem{definition}[theorem]{Definition}

\usepackage{comment}
\usepackage{tikz,pgfplots}
\usepackage{tkz-euclide}

\title{Blow-ups of minimal surfaces in the Heisenberg group}
\author{Yonghao Yu}
\thanks{NEW YORK UNIVERSITY, COURANT INSTITUTE OF MATHEMATICAL SCIENCES, 251 MERCER STREET, NEW YORK, NY 10012, USA. \\ 
Email: yy4004@nyu.edu}
\pgfplotsset{compat=1.18}
\begin{document}

\begin{abstract}
In this paper, we revise Monti's results \cite{monti_2015} on the blow-ups of H-perimeter minimizing sets in $\h^n$. Monti demonstrated that the Lipschitz approximation of the blow-up, after rescaling by the square root of the excess, converges to a limit function for $n \ge 2$. However, the partial differential equation he derived for this limit function $\varphi$ through contact variation is incorrect. Instead, the correct equation is that the horizontal Laplacian of the limit function $\varphi$ is independent of the coordinate $y_1$ and solves equation \ref{eq:thma} weakly. 
\end{abstract}
\maketitle
\section{introduction}

Let $E \subset \mathbb{H}^n$ be a $H$-perimeter minimizing set in some neighborhood of $0\in \h^n$ with $0 \in \partial^* E$, the $H$-reduced boundary of $E$. By rescaling the set $E$ with dilations $\delta_{1/r_h}$ as $r_h$ goes to $0$, Monti constructed a sequence of $H-$perimeter minimizing sets $E_h$ with horizontal excess $\eta_h^2$ approaching $0$.
Then, using the Lipschitz approximation theorem proved in \cite{monti_2013}, Monti obtained a sequence of intrinsic Lipschitz functions $\varphi_h:D\to \R$ that approximate the boundary of the rescaled sets $E_h$, where $D$ is some open subset of the vertical hyperplane $W=\{(x_1,y_1,\dots,x_n,y_n,t)\in \h^n: x_1=0\}$. As stated in Theorem \ref{thm: existence of limit function}, for $n\ge 2$, Monti showed that there exists a sequence of $h_i$ such that $(\varphi_h/\eta_h)$ weakly converges to $\varphi$ in $L^2(D)$. He claimed that the limit function $\varphi$ is independent of the first variable $y_1$. When $E$ is strongly perimeter minimizing (see Definition \ref{strong}), he claimed that $\varphi$ satisfied an equation involving the Kohn-Laplacian $\Delta_H=\sum_{i=2}^n X_i^2 + Y_i^2$. However, both claims are incorrect due to a calculation error. We correct his result in Theorem A.

We identify $\h^n$ with the set $\mathbb{C}^n \times \R$ by the coordinates $(z,t)$ where $z=(x_1+iy_1,\dots,x_n+iy_n)$. Let $X_i$, $Y_i$, $i=1, \dots, n$ denote the usual left-invariant vector fields defined on $\h^n$, and $Q_{r}$ be the homogeneous cube centered at $0$ with radius $r$, it is defined as 
$$
    Q_{r}= \{(z,t)\in \h^n: |x_i|<r,|y_i|<r,|t|<r^2, i=1,\dots,n\},
$$

\begin{theorema}\label{thma}
For $n\ge 2$, take any locally finite perimeter set $E\subset \h^n$. Suppose $0\in \partial^*E$ and the horizontal inner normal $\nu_E(0) = (1,0,\dots,0)\in \h^{n}$. Then:
\begin{enumerate}[$(i)$]
    \item If $E$ is $H$-perimeter minimizing in some neighborhood of $0\in \h^n$, the limit function $\varphi$ in Theorem \ref{thm: existence of limit function} solves the following equation weakly in $D_{1/4}$:
\begin{equation}\label{eq:thma}
    \frac{\partial}{\partial y_1}\Delta_0 \varphi = 0,
\end{equation}
where $\Delta_0= \frac{\partial^2}{\partial y^2_1} +\sum_{i=2}^n X_i^2 + Y_i^2$, $D_{1/4}=\{(z,t)\in Q_{1/4}: x_1=0\}$.
    \item If $E$ is strongly $H$-perimeter minimizing in some neighborhood of $0\in \h^n$, then the limit function $\varphi$ solves the following equation weakly in $D_{1/4}$:
\begin{equation}\label{eq:thma2}
    \Delta_0 \varphi=0 
\end{equation}
\end{enumerate}
\end{theorema}

\textbf{Acknowledgments:} This material is based upon work supported by the National Science Foundation under Award No. 2005609. The author would like to thank
Simone Verzellesi and Robert Young for their time and advice during the preparation of this paper.

\section{Preliminaries}

The Heisenberg group $\h^n$ is the set $\mathbb{C}^{n}\times \R$ equipped with the group product,
$$
 (z,t)* (z',t') = \left(z+z', t+t'-2\operatorname{Im}\langle z,z'\rangle\right),   
$$
where $z,z'\in \mathbb{C}^n$, $t,t'\in \R$, and $\langle z,z'\rangle=z_1 \bar{z}_1'+\dots +z_n \bar{z}_n'$. 

 The Lie algebra of $\h^n$ is spanned by the left-invariant vector fields,
$$
X_k = \frac{\partial}{\partial x_k} +2 y_k \frac{\partial}{\partial t},\quad Y_k= \frac{\partial}{\partial y_k} -2x_k \frac{\partial}{\partial t},\quad
T = \frac{\partial}{\partial t},
$$
where $z_k =x_k + iy_k$, $k=1,\dots, n$. The only non-trivial brackets are $[X_k,Y_k]=-4T$ for $k=1,\dots, n$. The vector fields $X_1,Y_1,\dots,X_n,Y_n$ are called the horizontal vector fields of $\h^n$. We define $H$ as the horizontal subbundle of $T\h^n$, where: $$H(p):= H_p = \text{Span}\{X_1(p),\dots, X_n(p),Y_1(p),\dots, Y_n(p)\}.$$

Let $\Omega\subset \h^n$ be an open set. Take any continuous function $V:\Omega \to \R^{2n}$. It can be identified with a horizontal vector field $V=\sum_{j=1}^n V_j X_j + V_{n+j} Y_j$, then the horizontal divergence of $V$ is
$$
\Div_H V = \sum_{j=1}^n X_j V_j + Y_j V_{n+j},
$$

In the Heisenberg group, one can define the $H$-perimeter similarly to the perimeter in Euclidean space. 

\begin{definition} ($H-$perimeter)
Let $E\subset \h^n$ be a measurable subset and $\Omega\subset \h^n$ be an open set. The \emph{$H$-perimeter} of $E$ in $\Omega$ is defined as:
$$
P_H (E; \Omega) = \sup\left\{\int_\Omega \chi_E(z,t) \Div_H V dz\,dt: V\in C_0^1(\Omega; \R^{2n}),\|V\|_{\infty}\le 1\right\},
$$
where $\chi_E(z,t)$ denotes the characteristic function of $E$. 
\end{definition}

We say that $E$ has \emph{finite $H$-perimeter} in $\Omega$ if $P_H (E; \Omega)<\infty$. Moreover, $E$ is said to have \emph{locally finite $H$-perimeter} in $\Omega$ if for any open set $A\subset\Omega$, $P_H(E;A)$ is finite. We denote $P_H (E; A)$ as $\mu_E(A)$ and view $\mu_E$ as a Radon measure $\mu_E$ on $\Omega$. The measure $\mu_E$ is called the $H$-perimeter measure of $E$. By the Riesz representation theorem, there exists a Borel function $\nu_E:\Omega \to \R^{2n}$ such that $|\nu_E|=1\mu_E$ a.e. and the following formula 
$$
\int_{\Omega}\left\langle V, \nu_E\right\rangle d \mu_E=-\int_{\Omega} \operatorname{div}_H V d z d t
$$ holds for any $V \in C_c^1\left(\Omega ; \mathbb{R}^{2 n}\right)$. We call the vector function $\nu_E$ the \emph{horizontal inner normal} of $E$ in $\Omega$. 

\begin{definition}
    A set $E$ with locally finite perimeter is considered to be \emph{$H$-perimeter minimizing} in a open set $U$ if:
\begin{equation}
    P_H(E,Q_r)\le P_H(F, Q_r),
\end{equation}
for any set $F\subset \h^n$ such that the symmetric difference $E\symdiff F$ is a compact subset of $U$.
\end{definition}

Let $Y_1$ be the vector field defined above, we define the closure of the cube $Q_r$ relative to the direction $Y_1$ as 
$$
\overline{Q}_r^{Y_1,+} = \left\{(z, t) \in \mathbb{H}^n : -r < y_1 \leq r\text{, } \left|x_1\right| < r\text{, }\left|t\right| < r^2, \text{ and } \left|x_i\right|, \left|y_i\right| < r \text{ for } i = 2, \ldots, n\right\}.
$$

\begin{definition}\label{strong}\cite{monti_2015} (strongly perimeter minimizing)
    A set $E$ with locally finite perimeter is considered to be  \emph{strongly $H$-perimeter minimizing in $Q_r$}, if for any $0<s\le r$, 
\begin{equation}
    P_H(E,Q_s)\le P_H(F, Q_s),
\end{equation}
for any set $F\subset \h^n$ such that $(E\symdiff F)\cap \bar{Q}_s$ is a compact subset of $\overline{Q}_s^{Y_1,+}$.
\end{definition}

\subsection{Excess and reduced boundary}

Let $E$ be a subset of $\mathbb{H}^n$ with locally finite $H$-perimeter, we call $0\in \h^n$ a point of the $H$-reduced boundary of $E$, denoted by $0\in \partial^*E$, if $\mu_E(B_r)>0$ for all $r>0$,
$$
\lim _{r \rightarrow 0} \frac{1}{\mu_E\left(B_r\right)} \int_{B_r} \nu_E\, d \mu_E=\nu_E(0),
$$
and $|\nu_E(0)|=1$. 

 Take any $p\in \h^n$, $r>0$, and $v\in S^{2n}$. The \emph{$v-$directional horizontal excess} of $E$ in $B_r(p)$ is 
\begin{equation}
\Exc(E,B_r(p),v)=\frac{1}{r^{2n+1}}\int_{B_r(p)}|\nu_E(p)-v|d\mu_E.
\end{equation}
The \emph{horizontal excess} of $E$ in $B_r(p)$ is the minimum of all possible directional horizontal excesses of $E$:
\begin{equation}
\Exc(E,B_r(p))=\min_{v\in S^{2n}}\Exc(E,B_r(p),v), 
\end{equation}

\subsection{Intrinsic Lipschitz graph}
 Recall $W=\{(z,t)\in \h^n, x_1=0\}=\R^{2n}$ is the vertical hyperplane. For any function $\varphi: W\to \R$, we define the \emph{intrinsic graph} of $\varphi$ along $X_1$ to be 
  \begin{equation}
        \gr(\varphi)=\{\left(z+\varphi(z,t)\mathit{e}_1, t+2 y_1 \varphi(z,t)\right): (z,t) \in W\},
  \end{equation}
  and the \emph{intrinsic epigraph} of $\varphi$ along $X_1$ to be 
\begin{equation}
    E_\varphi=\{\left(z+s\mathit{e}_1, t+2 y_1 s\right): (z,t) \in W, s>\varphi(z,t)\},
\end{equation}
where $z=(x_1,\dots,x_n,y_1,\dots,y_n)\in \mathbb{C}^n=\R^{2n}$ and $\mathit{e}_1=(1,0,\dots,0)\in \R^{2n}$.

Moreover, taking any point $w=(z,t)\in W$, we use the notation
\begin{equation}
    w*\varphi(w):=\left(z+\varphi(w)\mathit{e}_1, t+2 y_1 \varphi(w)\right)
\end{equation}
 Let $V=(\mathit{e}_1,0)\in \h^{n}$, for any $p\in \h^n$, set $\mathit{e}_1(p)=\langle p,V\rangle V\in \h^n$ to denote the projection of $p$ on to $\mathit{e}_1$, and let $\mathit{e}_1^{\perp}(p)$ be the point such that
\begin{equation}
    p = \mathit{e}_1^{\perp}(p)*\mathit{e}_1(p).
\end{equation}
Then the cone with vertex $0\in \h^{n}$, axis $\mathit{e}_1$ and aperture $\alpha \in (0,\infty]$ is the set
\begin{equation}
    C(0,v,\alpha) = \{p\in \h^n : \|v^{\perp}(p)\|_\infty < \alpha \|v(p)\|_\infty\}.
\end{equation}
Now, for a cone with vertex $p$ instead of $0$, we define the set as $C(p,v,\alpha)=p*C(0,v,\alpha)$.
\begin{definition}
 Let $D\subset W$ be an open set. A continuous function $\varphi: D\to \R$ is a \emph{$L$-intrinsic Lipschitz} function with $L \in[0, \infty)$, if for any $p \in \operatorname{gr}(\varphi)$ there holds
$$
\operatorname{gr}(\varphi) \cap C(p, v, 1 / L)=\emptyset.
$$
\end{definition}

The gradient of intrinsic Lipschitz function is called intrinsic gradient.
\begin{definition}
Let $D\subset W$ be an open set, 
for any function $\varphi\in \text{Lip}_{loc}(D)$, we define the \emph{intrinsic gradient} $\nabla^{\varphi} \varphi$ to be 
\begin{equation}
 \nabla^{\varphi} \varphi = 
\left(X_2 \varphi, \ldots, X_n \varphi, \mathfrak{B} \varphi, Y_2 \varphi, \ldots, Y_n \varphi\right),  
\end{equation}
where $\mathfrak{B}$ is the Burgers' operator,
\begin{equation}\label{Burger}
    \mathfrak{B}\varphi= \frac{\partial \varphi}{\partial y_1}-4\varphi \frac{\partial \varphi}{\partial t}.
\end{equation}

If $\varphi \in C(D)$ is a continuous function, we say that the intrinsic gradient $\nabla^{\varphi} \varphi$ exist in the sense of distributions if $X_i\varphi$, $\mathfrak{B}\varphi$, $Y_i\varphi$, $i=2,\dots, n$ exists in the sense of distributions. Then $\nabla^{\varphi} \varphi \subset L_{loc}^{\infty}(D; \R^{2n-1}).$  
\end{definition}

In \cite{monti_2013}, Monti showed that the boundary of a set of minimizing H parameters $E$ can be approximated with a $L-$ Lipschitz graph. Let $\mathscr{S}^{2n+1}$ denote the $(2n+1)$ dimensional spherical Hausdorff metric associated with the Carnot-Carathedory distance. Then there is an intrinsic Lipschitz function $\varphi$ such that the measure of the symmetric difference of $\gr(\varphi)$ and $\partial E$ is bounded by the excess, as the theorem below shows.
\begin{theorem}[Lipschitz approximation]\label{thm:LA} \cite[Theorem~1.1]{monti_2013} Let $n\ge 2$. For $L>0$, there exists some constant $k>1$ such that for any H-perimeter minimizing set $E$ in $B_{kr}$, with $0\in \partial E$ and $r>0$, there exists an $L-$intrinsic Lipschitz function $\varphi: W\to \R$ such that
\begin{equation}
\mathscr{S}^{2n+1}\left((\operatorname{gr}(\varphi) \symdiff \partial E) \cap B_r\right) \leq c(L, n)(k r)^{Q-1} \operatorname{Exc}\left(E, B_{k r},X_1\right),
\end{equation}
where $c$ is some positive constant which depends on $L$ and $n$. 
\end{theorem}

\subsection{Contact flow and First Variation}

Taking any bounded open set $\Omega\in \h^n$, a \emph{contact map} defined on $\Omega$ is a diffeomorphism $\Psi$ such that the differential map $\Psi_*$ preserves the horizontal subbundle. That is, for any $p\in \Omega$, $\Psi_* (H_p)\subset H_{\Psi(p)}$. Then a one-parameter flow $(\Psi_s)_{s\in\R}$ of $\h^n$ is a \emph{contact flow} if each $\Psi_s$ is a contact map. For more information on contact flows, see \cite{koranyi_reimann_1995}.

Take any generating function $\psi\in C^\infty (\h^n)$, let $V_{\psi}$ be the vector field in $\h^n$ of the form
\begin{equation}\label{contactvectorfield}
   V_{\psi} = \sum_{j=1}^n (Y_j \psi) X_j - (X_j \psi)Y_j - 4\psi T.
\end{equation} 
Then there is a corresponding contact flow $\Psi: [-\delta, \delta]\times \Omega\to \h^n$, such that for any $s\in [-\delta, \delta]$, and $p\in \Omega$, the following relation holds:
\begin{align}
    \Psi'(s,p)&= V_{\psi}(\Psi(s,p)),\\
    \Psi(0,p)&=p.
\end{align}
we call the flow generated by some generating function $\psi\in C^\infty (\h^n)$.

In \cite{monti_2014}, Monti showed the following first variation formula for a contact flow $\Psi$.
\begin{theorem}\cite[Theorem~3.18]{monti_2014} Let $\Omega$ be some bounded open set in $\h^n$, and $\Psi: [-\delta, \delta]\times \Omega\to \h^n$ be a contact flow generated by some smooth function $\psi\in C^{\infty}(\h^n)$. Then there exists some positive constant $C$ depending on $\psi$ and $\Omega$ such that for any set $E\subset \h^n$ with a finite perimeter in $\Omega$, we have 
\begin{equation}\label{FV}
\left|P_H\left(\Psi_s(E), \Psi_s(\Omega)\right)-P_H(E, \Omega)+s \int_{\Omega}\left\{4(n+1) T \psi+L_\psi\left(\nu_E\right)\right\} d \mu_E\right| \leq C P_H(E, \Omega) s^2,
\end{equation}
for any $s\in [-\delta, \delta]$, where $L_\psi: H_p\to \R$ is the real quadratic form,
\begin{equation}\label{eq:L}
L_\psi\left(\sum_{j=1}^n x_j X_j+y_j Y_j\right)=\sum_{i, j=1}^n x_i x_j X_j Y_i \psi+x_j y_i\left(Y_i Y_j \psi-X_j X_i \psi\right)-y_i y_j Y_j X_i \psi.
\end{equation}
\end{theorem}

\section{Proof of Theorem A}

\quad We need to recall an approximation result before the proof of theorem A. Let $E\subset \h^n$ be a $H$-perimeter minimizing set in some neighborhood of $0$, with $0\in \partial^*E$, $\nu_E(0) = X_1$. Since $0\in \partial^*E$, there exists a sequence of real numbers $r_h\to 0^+$ such that
$$\text{Exc}(E, Q_{r_h})<\frac{1}{h}.$$

    For any real number $\lambda>0$, let $\delta_\lambda(z,t)=(\lambda z,\lambda^2 t)$ be the dilation defined on $\h^n$. Then the rescaled sets $E_h=\delta_{1/r_h}(E)$ satisfy the following properties:
\begin{enumerate}[(i)]
    \item $0\in \partial^* E_h$, and $\nu_{E_h} (0)=X_1$ .
    \item Each set $E_h$ is $H$-perimeter minimizing.
    \item $\text{Exc}(E_h,Q_1)<\frac{1}{h}$ since excess is dilation invariant.
\end{enumerate}

Let $\eta_h = \sqrt{\Exc(E_h,Q_1)}$ denote the square root of the excess of $E_h$. Pick some small number $\sigma$ such that $0<\sigma<\frac{1}{k}$, where $k$ is the geometric constant defined in Theorem \ref{thm:LA}. Then by Theorem \ref{thm:LA}, there exists a $L-$intrinsic Lipschitz function $\varphi_h : W\to \R$ such that 
\begin{equation}\label{eq: 2.22}
\begin{aligned}
   \mathscr{S}^{Q-1}\left(\left(\operatorname{gr}\left(\varphi_h\right) \symdiff \partial E_h\right) \cap B_\sigma\right)&\le \mathscr{S}^{Q-1}\left(\left(\operatorname{gr}\left(\varphi_h\right) \symdiff \partial E_h\right) \cap B_{1/k}\right)\\
   &\leq c(L, n,\sigma)(\frac{1}{k}\cdot k)^{n+1} \operatorname{Exc}\left(E_h, B_1\right)=c_0 \eta_h^2, 
\end{aligned}
\end{equation}
where $c_0=c(L,n,\sigma)$.

Using this inequality and a Poincaré-type inequality proved in \cite{cmpsc_Poincare_2016}, Monti showed that there exists a subsequence of $(\varphi_h/\eta_h)_{h\in\N}$ that weakly converges to some function $\varphi$ in the $L^2$ sense. 
\begin{theorem}\label{thm: existence of limit function} \cite[Theorem~ 2.5]{monti_2015}
 Assume $n\ge 2$, following from the construction above, let $\varphi_h$ be the L-intrinsic function associated with $E_h$. Then there exists an open neighborhood $D\subset W$ of $0$, real constants $\overline{\varphi_h}$, and a selection of indices $k\to h_k$ such that as $k\to \infty$, $\frac{\varphi_{h_k}-\overline{\varphi_{h_k}}}{\eta_h}$  weakly converge to some function $\varphi \in W^{1,2}_{H}(D)$. Moreover, the intrinsic gradient $\nabla^{\varphi_{h_k}} \varphi_{h_k}$ also converges:
\begin{equation}
\frac{\nabla^{\varphi_{h_k}} \varphi_{h_k}}{\eta_{h_k}} \rightarrow \nabla_H \varphi \text { weakly in } L^2\left(D ; \mathbb{R}^{2 n-1}\right),    
\end{equation}
 
where 
\begin{align}
    \nabla_H \varphi = \left(X_2\varphi,\dots, X_n\varphi,\frac{\partial \varphi}{\partial y_1},Y_2\varphi,\dots, Y_n \varphi\right)
\end{align}
\end{theorem}
\begin{remark}
     In the proof of theorem 2.1, Monti also showed that
\begin{equation}\label{22}
    \varphi_{h}\to 0\text{ and }\overline{\varphi_h}\to 0\text{ strongly in }L^2(D)
\end{equation}
and the following estimate of intrinsic gradient
\begin{equation}
    \int_{D_1} |\nabla^{\varphi_h}\varphi_h|^2\le c_0 \eta_h^2, \label{eta_h^2}
\end{equation}
where $c_0$ is the constant in (\ref{eq: 2.22}). 
\end{remark}

Without loss of generality, we can assume that $D_1=\left\{(z,t)\in Q_1:x_1=0\right\}\subset D$, then the limit function $\varphi$ is defined on the whole $D_1$.

\noindent\textbf{Proof of Theorem A:}
The proof is a revised version of Monti's proof of Theorem $3.2$ in \cite{monti_2015}, where he used contact flow and first variation formula to obtain the final result. He assumed that the generating function $\psi\in C^{\infty}(\h^n)$ for the contact flow $\Psi:[-\delta,\delta]\times D_1 \to \h^n$ is of the form
$$
\psi = \alpha + x_1\beta + \frac{1}{2}x_1^2\gamma,  
$$
where $\alpha, \beta, \gamma$ are smooth functions in $\h^n$ such that $$X_1\alpha=X_1 \beta = X_1 \gamma = 0 \text{ in } Q_{1/2}.$$ 

He further assumes that $\beta, \gamma$ are compactly supported in $Q_{1/2}$. We can identify two cases in Theorem A just as Monti did in $(3.45)$ and $(3.46)$: 

\begin{enumerate}[$(i)$]
    \item If $E$ is $H$-perimeter minimizing in $Q_1$, we assume:
    \begin{equation}
        \alpha \in C_c^{\infty}\left(Q_{1 / 2}\right),
    \end{equation}
    Then by (\ref{contactvectorfield}), the contact vector field $V_\psi$ vanishes outside $Q_{1/2}$. We have $E_h \triangle \Psi_s(E_h) \subset Q_1$ for any $s\in (0,\delta]$. Hence by the definition of $H$-perimeter minimizing, $P_H(E_h,Q_1)\le P_H(\Psi_s(E_h),Q_1)$ .     
    \item If $E$ is strongly $H$-perimeter minimizing in $Q_1$, we will first define $\alpha_0$ on the vertical hyperplane $W$:
\begin{equation}
    \alpha_0\left(y_1, z_2, \ldots, z_n, t\right)=\int_0^{y_1} \vartheta_0\left(s, z_2, \ldots, z_n, t\right) d s, \quad y_1 \in \mathbb{R}\text{, }z_2,\dots,z_n\in\mathbb{C},
\end{equation}
where $\vartheta_0 \in C_c^{\infty}\left(D_{1 / 2}\right)$ and $X_1 \vartheta_0=0$.
Moreover, letting $\Pi: \h^n\to W$ denote the nonlinear projection along the cosets of $<X_1>$ given by $\Pi(z_1,\dots,z_n,t)=(0+iy_1,z_2,\dots,z_n,t-2x_1y_1)$, we define $\alpha(z_1,\dots,z_n,t):=\alpha_0\circ\Pi(z_1,\dots,z_n,t)$. Then $X_1 \alpha =0$ in $\h^n$, $\alpha$ is supported in the $y_1$ cylinder $C_{y_1}=\{(z,t)\in \h^n:|x_1|<\frac{1}{2},|z_i|<\frac{1}{2},|t|<\frac{1}{4}, i=2,\dots,n\}$, and remains constant for $|y_1|\ge \frac{1}{2}$. Then $V_\psi$ is supported in $C_{y_1}$ and equals $-4\alpha T$ for $|y_1|\ge \frac{1}{2}$. Then we have $E_h \triangle \Psi_s(E_h)\cap \overline{Q_1} \subset \overline{Q}_1^{Y_1,+}$ for any $s\in(0,\delta]$. Hence, by the definition of strongly $H$-perimeter minimizing, $P_H(E_h,Q_1)\le P_H(\Psi_s(E_h),Q_1)$.
\end{enumerate}

With abuse of notation, we use $\psi,\alpha,\beta,\gamma$ to denote its restriction on the plane $\{x_1=0\}$. We also use the notations:
$$f_t=\frac{\partial f}{\partial t},$$ 
$$f_{y_1}=\frac{\partial f}{\partial y_1}$$
for any smooth functions $f$ defined on $D$.

Then Monti applied the first variation formula (\ref{FV}) to each rescaled set $E_h$. By the minimality condition $P_H(E_h,Q_1)\le P_H(\Psi_s(E_h),Q_1)$ and the weak convergence of $\varphi_h/\eta_h$, Monti concluded $(3.49)$ in \cite{monti_2015}:
\begin{equation}\label{eq: first_variational_0}
    \Delta P_h:=\lim_{h\to \infty} \frac{1}{\eta_h}\int_D \left\{4(n+1)T\psi(w*\varphi_h(w))+L_{\psi}\left(\nu_{E_{\varphi_h}}(w*\varphi_h(w))\right)\right\}dw = 0,
\end{equation}
where $D$ denotes the unit disk $D_1$ on the vertical plane $W$, $L_{\psi}$ is the quadratic form (\ref{eq:L}) associated with the first variation, $E_{\varphi_h}$ is the intrinsic epigraph of $\varphi_h$ and $\nu_{E_{\varphi_h}}$ is the horizontal inner normal of $E_{\varphi_h}$.

As shown from $(3.51)$ to $(3.53)$ in \cite{monti_2015}, Monti computed the first half of $\Delta P_h$ as follows:
\begin{equation}
    \lim_{h\to \infty} \frac{1}{\eta_h}\int_D 4(n+1)T\psi(w*\varphi_h(w))dw=\int_ D 4(n+1)\beta_t \varphi.
\label{eq:1}    
\end{equation}
It remains to calculate the second half of $\Delta P_h$
\begin{equation}
    \lim_{h\to \infty} \frac{1}{\eta_h}\int_D L_{\psi} \left(\nu_{E_{\varphi_h}}(w*\varphi_h(w))\right)dw,
    \label{eq:2}
\end{equation}

The error occurred when Monti was expanding the form $L_{\psi}(\nu_{E_{\varphi_h}})$, on the third line of $(3.56)$ in \cite{monti_2015}. The term $x_1 X_j \gamma$ should be $\gamma$ instead of $x_1 X_1 \gamma$ when $j=1$. This caused the final integral $(3.62)$ in \cite{monti_2015} to miss a term $\gamma \varphi_{y_1}$. 

By Theorem 2.4 in \cite{monti_2015}, the horizontal normal $\nu_{E_{\varphi_h}}=(\nu_{X_1},\dots, \nu_{X_n},\nu_{Y_1},\dots \nu_{Y_n})$ is of the form:
\begin{equation}\label{31}
    \nu_{X_1}= \frac{1}{\sqrt{1+|\nabla^{\varphi_h}\varphi_h}|^2},  \nu_{Y_1}= -\frac{\mathfrak{B}\varphi_h}{\sqrt{1+|\nabla^{\varphi_h}\varphi_h}|^2},
\end{equation}
\begin{equation}\label{32}
   \nu_{X_i}= -\frac{X_i \varphi_h}{\sqrt{1+|\nabla^{\varphi_h}\varphi_h}|^2},  \nu_{Y_i}= -\frac{Y_i \varphi_h}{\sqrt{1+|\nabla^{\varphi_h}\varphi_h}|^2},  
\end{equation}
for $2\le i\le n$, where $\mathfrak{B}$ is the Burgers operator.

 Then Monti rewrote $L_{\psi}$ as terms that contain $\nu_{X_1}$ plus some quadratic form $K_{\psi}(\nu_{E_{\varphi_h}})$
\begin{equation}\label{eq:Lpsi}
    \begin{aligned}
    L_{\psi}\left(\nu_{E_{\varphi_h}}\right) = \nu_{X_1}^2 X_1 Y_1 \psi +\sum_{i=2}^n \left(\nu_{X_i} \nu_{X_1} X_1 Y_i \psi + \nu_{X_1}\nu_{X_i}X_i Y_1 \psi\right)\\
    +\sum_{i=1}^n \nu_{X_1}\nu_{Y_i}(Y_i Y_1 \psi - X_1 X_i \psi)+ K_\psi\left(\nu_{E_{\varphi_h}}\right).
    \end{aligned} 
\end{equation}

    Then he computed the derivatives of $\psi$ in $L_\psi$ that associate with $\nu_{X_1}$:
\begin{equation}
\begin{aligned}
X_1 Y_1 \psi &= Y_1 X_1 \psi - 4T \psi\\
& = Y_1 \beta + x_1 Y_1 \gamma -4 \left(\alpha_t + x_1 \beta_t +\frac{1}{2}x_1 ^2 \gamma_t\right),\\
X_1 Y_i \psi & = Y_i X_1 \psi = Y_i \beta + x_1 Y_i \gamma \quad(i\ge 2),\\
Y_1 Y_i \psi &= Y_i Y_1 \alpha + x_1 Y_i Y_1 \beta + \frac{1}{2}x_1^2 Y_i Y_i \gamma  \quad(i \ge 1),\\
X_1X_1 \psi &=X_1(\beta + x_1 \gamma)=\gamma,\\
X_1 X_i \psi &= X_i X_1 \psi = X_i \beta + x_1 X_i \gamma \quad (i\ge 2), \\
X_i Y_1 \psi &= X_iY_1 \alpha + x_1 X_i Y_1\beta + \frac{1}{2}x_1^2 X_i Y_1 \gamma  \quad(i\ge 2).
\end{aligned}   
\end{equation}
    Inserting the derivatives into (\ref{eq:Lpsi}), we obtain the corrected version of $(3.56)$, 
\begin{equation}
\begin{aligned}
L_\psi\left(\nu_{E_{\varphi_h}}\right)=L_1(h)+L_2(h)+L_3(h)+L_4(h)+K_\psi\left(\nu_{E_{\varphi_h}}\right).\label{48}
\end{aligned}
\end{equation}

where we define
\begin{align}
    L_1(h)&:=\left\{Y_1 \beta+x_1 Y_1 \gamma-4\left(\alpha_t+x_1 \beta_t+\frac{1}{2} x_1^2 \gamma_t\right)\right\} \nu_{X_1}^2,\\
    L_2(h)&:=\sum_{j=2}^n\left\{Y_1 X_j \alpha+x_1 Y_1 X_j \beta+\frac{1}{2} x_1^2 Y_1 X_j \gamma+Y_j \beta+x_1 Y_j \gamma\right\} \nu_{X_1} \nu_{X_j}, \\
    L_3(h)&:=\sum_{j=2}^n\left\{Y_j Y_1 \alpha+x_1 Y_j Y_1 \beta+\frac{1}{2} x_1^2 Y_j Y_1 \gamma-X_j \beta-x_1 X_j \gamma\right\} \nu_{X_1} \nu_{Y_j},\\
    L_4(h)&:=\left(Y_1Y_1 \alpha+x_1 Y_1Y_1 \beta+\frac{1}{2} x_1^2 Y_1Y_1 \gamma-\gamma\right) \nu_{X_1} \nu_{Y_1}.
\end{align}

Since $\varphi_h$ is $L$-intrinsic Lipschitz, we can assume that its intrinsic gradient $\nabla^{\varphi_h}\varphi_h$ must be bounded everywhere. Hence there exists some large constant $A>0$ such that 
\begin{equation}
\left|K_{\varphi}(\nu_{E_{\varphi_h}})\right|\le A\left|\nabla^{\varphi_h}\varphi_h\right|^2 .  
\end{equation}
By (\ref{eta_h^2}), we have $ \int_D|\nabla^{\varphi_h}\varphi_h|^2 \le c_0\eta_h^2 $ for some positive constant $c_0$, then we obtain
\begin{equation}
     \lim_{h\to \infty} \frac{1}{\eta_h}\int_D \left|K_{\psi} \left(\nu_{E_{\varphi_h}}(w*\varphi_h(w))\right)\right|dw\le \lim_{h\to \infty}\frac{1}{\eta_h}\int_D A\left|\nabla^{\varphi_h}\varphi_h\right|^2 dw \le   \lim_{h\to \infty} \frac{A}{\eta_h} c_0\eta_h^2 =0.
\end{equation}

Next, we compute the limit of the integral of $L_1(h)$, $L_2(h)$, $L_3(h)$ and $L_4(h)$. Since $X_1 \beta = 0$, $X_1 Y_1 \beta = Y_1 X_1 \beta -4T \beta = -4T \beta$, and noticing that $x_1=\varphi_h$ we have 
\begin{equation}
Y_1\beta(w*\varphi_h(w))=Y_1\beta(w)-4x_1T\beta(w)=\beta_{y_1}(w)-4\varphi_h \beta_t(w).
\label{eq:beta}    
\end{equation}

Similarly, one can show that 
\begin{equation}
Y_1\gamma(w*\varphi_h(w))=\gamma_{y_1}-4\varphi_h \gamma_t , 
\label{eq:gamma}
\end{equation}
\begin{equation}
Y_1\alpha(w*\varphi_h(w))=\alpha_{y_1}-4\varphi_h \alpha_t  .
\label{eq:alpha}
\end{equation}

Let $\overline{\varphi_h}$ be the real constant in Theorem \ref{thm: existence of limit function}, by (\ref{22}) we have
\begin{equation}\label{x^2}
\lim_{h\to \infty}\int_D \frac{\varphi_h^2}{\eta_h}f=\lim_{h\to \infty}\int_D \frac{\varphi_h^2-\overline{\varphi_h}^2}{\eta_h}f=\lim_{h\to \infty}\int_D \frac{\varphi_h-\overline{\varphi_h}}{\eta_h}(\varphi_h+\overline{\varphi_h})f=0    
\end{equation}
for any smooth function $f\in C^\infty (D)$.

Moreover, since $\alpha$ is compactly supported in $|t|\le \frac{1}{2}$, $\beta$ and $\gamma$ is compactly supported in $Q_{1/2}$, we have
\begin{equation}\label{beta}
    \int_D \beta_{y_1}=\int_D \gamma_{y_1}=\int_D \alpha_t=\int_D \beta_t=0
\end{equation}

 Noticing that $\nu_{X_1}=1$ as $h$ goes to infinity by (\ref{31}). Then by theorem \ref{thm: existence of limit function},  (\ref{eq:beta}), (\ref{x^2}), and (\ref{beta}), the limit of the integral of $L_1(h)$ becomes
\begin{equation}
\begin{aligned}
 \lim_{h\to \infty}\frac{1}{\eta_h}\int_D  L_1(h) dw  &=\lim_{h\to \infty} \frac{1}{\eta_h}\int_D \left((\beta_{y_1}-4\varphi_h \beta_t) + \varphi_h(\gamma_{y_1}-4\varphi_h \gamma_t ) -4 \left(\alpha_t + \varphi_h \beta_t +\frac{1}{2}\varphi_h ^2 \gamma_t\right)\right)\nu_{X_1}^2dw\\
&= \lim_{h\to \infty} \int_D (-8\beta_t +\gamma_{y_1})\frac{\varphi_h}{\eta_h}dw\\
&= \lim_{h\to \infty} \int_D (-8\beta_t +\gamma_{y_1})\frac{\varphi_h-\overline{\varphi_h}}{\eta_h}dw\\
&=\int_D (\gamma_{y_1}-8\beta_t)\varphi dw.
\end{aligned}
\label{eq:a}
\end{equation}

Similarly, using the fact that $\lim_{h\to \infty}\varphi_h = 0$, theorem \ref{thm: existence of limit function}, (\ref{eq:alpha}) and (\ref{32}), the limit of integral of $L_2(h)$ becomes,
\begin{equation} 
\begin{aligned}
   \label{eq:b}
    \lim_{h\to \infty}\frac{1}{\eta_h}\int_D  L_2(h) dw  &= \lim_{h\to \infty}\frac{1}{\eta_h}\int_D \sum_{i=2}^n(Y_i \beta + \varphi_h Y_i \gamma+X_iY_1 \alpha + \varphi_h X_i Y_1\beta + \frac{1}{2}\varphi_h^2 X_i Y_1 \gamma ) {\nu_{X_1}\nu_{X_i}} dw\\
   &= \lim_{h\to \infty} \int_D \sum_{i=2}^n (Y_i \beta +X_i Y_1 \alpha) \frac{\nu_{X_1}\nu_{X_i}}{\eta_h} dw\\
    &= - \lim_{h\to \infty} \int_D \sum_{i=2}^n (Y_i \beta +X_i (\alpha_{y_1}-\varphi_h \alpha_t) \frac{X_i \varphi_h}{\eta_h} dw\\
   & =  -\int_D \sum_{i=2}^n (Y_i \beta +X_i  \alpha_{y_1})X_i\varphi dw,
   \end{aligned}
\end{equation}
Furthermore, one can use the same logic to compute the limit of the integral of $L_3(h)$
\begin{equation} 
\begin{aligned}
   \label{eq:c2}
      \lim_{h\to \infty}\frac{1}{\eta_h}\int_D  L_3(h) dw  &=\lim_{h\to \infty}\frac{1}{\eta_h} \int_D \sum_{i=2}^n(Y_i Y_1 \alpha + \varphi_h Y_i Y_1 \beta + \frac{1}{2}\varphi_h^2 Y_i Y_i \gamma -(X_i \beta + \varphi_h X_i\gamma ) \nu_{X_1}\nu_{X_i} dw\\
   &= \lim_{h\to \infty} \int_D \sum_{i=2}^n (Y_i Y_1\alpha - X_i \beta) \frac{\nu_{X_1}\nu_{Y_i}}{\eta_h} dw\\ 
   &= - \lim_{h\to \infty} \int_D \sum_{i=2}^n (Y_i (\alpha_{y_1}-4\varphi_h \alpha_t) - X_i \beta) \frac{Y_i\varphi_h}{\eta_h} dw\\
   &= - \lim_{h\to \infty} \int_D \sum_{i=2}^n (Y_i \alpha_{y_1} - X_i \beta) \frac{Y_i\varphi_h}{\eta_h} dw\\
   & =  \int_D \sum_{i=2}^n (-Y_i \alpha_{y_i} +X_i \beta)Y_i\varphi dw.
   \end{aligned}
\end{equation}

Finally, we compute the limit of the integral of $L_4(h)$, which is the special case of Monti's computation in (3.61) in \cite{monti_2015} where $j=1$:
\begin{equation} 
\begin{aligned}
   \label{eq:c1}
  \lim_{h\to \infty}\frac{1}{\eta_h}\int_D  L_4(h) dw   &=\lim_{h\to \infty} \frac{1}{\eta_h}\int_D (Y_1 Y_1 \psi - X_1 X_1 \psi)\nu_{X_1} \nu_{Y_1} dw\\
   & = \lim_{h\to \infty}\int_D \left(Y_1^2 \alpha +\varphi_h Y_1^2 \beta + \frac{1}{2}\varphi_h^2 Y_1^2\gamma -\gamma\right)\frac{\nu_{X_1}\nu_{Y_1}}{\eta_h} dw\\
   &= -\lim_{h\to \infty}\int_D (Y_1^2\alpha - \gamma)\frac{\mathscr{B} \varphi_h}{\eta_h} dw\\
   & = -\int_D (Y_1^2 \alpha -\gamma) \varphi_{y_1} dw.\\
   &= \int_D -\varphi_{y_1} Y_1^2 \alpha  +\gamma \varphi_{y_1} dw.
   \end{aligned}
\end{equation}

 We obtain an additional $\gamma \varphi_{y_1}$ term in the integral (\ref{eq:c1}) compared to Monti's original integral $(3.61)$. 

Combining together (\ref{eq:1}) and  (\ref{eq:a}) to (\ref{eq:c1}),  \eqref{eq: first_variational_0} becomes:
\begin{equation}\label{eq:f}
\begin{aligned}
 \int_D (4(n-1)\beta_t+\gamma_{y_1} )\varphi - \varphi_{y_1} Y_1^2 \alpha+\gamma\varphi_{y_1}-\sum_{i=2}^n [(X_i\alpha_{y_1}+Y_i\beta)X_i \varphi
 +(Y_i\alpha_{y_1}-X_i \beta)Y_i \varphi]dw = 0.
\end{aligned}    
\end{equation}

 Note that this is $(3.62)$ in \cite{monti_2015} with the extra $\gamma \varphi_{y_1}$ term.

Setting $\alpha = \beta = 0$, we have 
\begin{equation}
\begin{aligned}
0 &= \int_D \gamma_{y_1} \varphi +\gamma\varphi_{y_1} dw\\
 & = \int_D (\gamma\varphi)_{y_1} dw.
\end{aligned}    
\end{equation}
which gives empty information since $\gamma\varphi$ is a compactly supported function on $D$. Therefore, $\int_D (\gamma\varphi)_{y_1}$ automatically equals zero. We no longer obtain claim i) in Theorem 3.2 in \cite{monti_2015} because we no longer obtain the formula:
$
0=\int_D \gamma_{y_1} \varphi d w=-\int_D \gamma \varphi_{y_1} dw
$.

Letting $\beta=\gamma =0$, using integration by parts, we obtain
\begin{equation}
\begin{aligned}
0 & = -\int_D \frac{\partial^2 \alpha}{\partial y_1^2} \varphi_{y_1}+ \sum_{i=2}^n (X_i \alpha_{y_1} X_i \varphi +Y_i \alpha_{y_1} Y_i \varphi) dw\\
& = \int_D \alpha_{y_1} \left\{\frac{\partial^2 \varphi}{\partial y_1^2}+\sum_{i=2}^n (X_i^2 \varphi + Y_i^2 \varphi)\right\} dw\\
& = \int_D \alpha_{y_1} \Delta_0 \varphi dw\\
& = -\int_D \alpha \frac{\partial}{\partial y_1}\Delta_0 \varphi dw.
\end{aligned}    
\end{equation}
If $E$ is perimeter minimizing, this holds for any compactly supported function $\alpha$,and we obtain the following differential equation:
$$
\frac{\partial}{\partial y_1}\Delta_0 \varphi = 0.
$$

If $E$ is strongly perimeter minimizing, $\alpha_{y_1}=\vartheta$ for any test function $\vartheta\in C^{\infty}_c(D_{1/2})$, then we have
$$
0= \int_D \vartheta \Delta_0 \varphi dw.
$$
Hence, $\varphi\in W_H^{1,2}(D)$ solves the partial differential equation $\Delta_0 \varphi=0$ in the weak sense.

\bibliographystyle{amsplain}
\bibliography{reference}
\end{document}